\documentclass[12pt]{article}
\pdfoutput=1
\usepackage[margin=1.25in]{geometry}
\usepackage{amscd}
\usepackage{amssymb,amsmath}
\usepackage{amsthm}
\usepackage{hyperref,xcolor}
\usepackage{graphicx}
\usepackage{cases}
\usepackage{indentfirst}
\usepackage{mathrsfs}
\usepackage{textcomp}
\usepackage{enumerate}
\usepackage[marginal]{footmisc}
\hypersetup{
	colorlinks,%
	citecolor=blue,%
	filecolor=blue,%
	linkcolor=blue,%
	urlcolor=black
}

\usepackage{verbatim}
\usepackage{sectsty}

\makeatletter

\newtheorem{theorem}{Theorem}[section]

\newtheorem{lemma}[theorem]{Lemma}
\newtheorem{remark}[theorem]{Remark}
\newtheorem{proposition}[theorem]{Proposition}
\newtheorem{corollary}[theorem]{Corollary}

\numberwithin{equation}{section}

\newdimen\bibspace
\renewenvironment{thebibliography}[1]{%
	\section*{\refname 
		\@mkboth{\MakeUppercase\refname}{\MakeUppercase\refname}}%
	\list{\@biblabel{\@arabic\c@enumiv}}%
	{\settowidth\labelwidth{\@biblabel{#1}}%
		\leftmargin\labelwidth
		\advance\leftmargin\labelsep
		\itemsep\bibspace
		\parsep\z@skip     %
		\@openbib@code
		\usecounter{enumiv}%
		\let\p@enumiv\@empty
		\renewcommand\theenumiv{\@arabic\c@enumiv}}%
	\sloppy\clubpenalty4000\widowpenalty4000%
	\sfcode`\.\@m}
{\def\@noitemerr
	{\@latex@warning{Empty `thebibliography' environment}}%
	\endlist}

\makeatother
\makeatletter

\def\XXint#1#2#3{{\setbox0=\hbox{$#1{#2#3}{\int}$}
		\vcenter{\hbox{$#2#3$}}\kern-.5\wd0}}

\newcommand{\be}{\begin{equation}}      \newcommand{\ee}{\end{equation}}

\begin{document}	
	\title{\bf\Large Existence and asymptotic behavior of entire large solutions for Hessian equations 
		\footnotetext{\hspace{-0.35cm}
			Xiang Li
			\endgraf 202031130022@mail.bnu.edu.cn
			\vspace{0.25cm}
			\endgraf Corresponding author: Jiguang Bao
			\endgraf jgbao@bnu.edu.cn
						\vspace{0.25cm}
			\endgraf All authors are supported in part by  the Beijing Natural Science Foundation (1222017).
	}}
	\vspace{0.25cm}
	\author{Xiang Li,\ \  Jiguang Bao\footnote{School of Mathematical Sciences, Beijing Normal University,
			 Beijing 100875, China}}
	\date{}
	\maketitle
\vspace{-0.8cm}
\noindent{\bf Abstract}		
	 In this paper, we give some existence  and nonexistence results for nonradial entire large solutions of the Hessian equation $S_k\left(D^2 u\right)=b(x) u^\gamma$ in the sublinear case $0<\gamma<k$. The exact asymptotic behavior of  large solutions at infinity is also studied when  $b(x)$ is the  oscillation of a radial function   $|x|^{-l}$ at infinity for $l\leq k-1$.

\noindent{\bf Keywords}  Hessian equation $\cdot$ sublinear $\cdot$ existence $\cdot$ asymptotic behavior $\cdot$ large solution 

\noindent{\bf Mathematics Subject Classification} $35{\rm J}60 \cdot 35{\rm A}01 \cdot 35{\rm B}40$	
	\section{Introduction and the main results}
	In this paper, we analyze the existence and asymptotic behavior  of entire  solutions of the following  blow-up problem
	\begin{equation}\label{y101}
	\begin{gathered}
	S_k\left(D^2 u\right):=\sigma_k\left(\lambda\left(D^2 u\right)\right)=b(x) u^\gamma, ~x \in \mathbb{R}^n, 
	\end{gathered}
	\end{equation}
	\begin{equation}\label{y113}
		u(x) \rightarrow \infty,~|x| \rightarrow \infty,
	\end{equation}
	where $b(x)$ is a positive continuous function in $\mathbb{R}^n$, $n\geq 3$, $\gamma \in(0,k)$,  and 
 $\lambda=\left(\lambda_1, \lambda_2, \cdots, \lambda_n\right)$ are the eigenvalues of the symmetric matrix $D^2 u$ and
	\begin{equation*}
	\sigma_k(\lambda):=\sum_{1 \leq i_1<\cdots<i_k \leq n} \lambda_{i_1}\cdots \lambda_{i_k}, ~ k=1,2,\cdots,n,
	\end{equation*}
	is the $k$-th elementary symmetric function.
	  Such a solution is called a large solution.

	To work in the realm of elliptic operators, we have to restrict the class of functions.  $u(x)$ is called a  $k$-convex function in $\Omega$ if it belongs to
		\begin{equation*}
	\Phi^k(\Omega):=\left\{u\left(x\right) \in C^2(\Omega): \lambda\left(D^2 u\left(x\right)\right) \in \Gamma_k, ~ x \in \Omega\right\},
	\end{equation*}
	where $\Omega$ is a domain in $\mathbb{R}^n$ and
	\begin{equation*}
	\Gamma_k:=\left\{\lambda \in \mathbb{R}^n: \sigma_l(\lambda)>0,1 \leq l \leq k\right\}.
	\end{equation*}
The $k$-Hessian operator $\sigma_k\left(\lambda\left(D^2 u \right)\right)$, introduced by Trudinger and Wang in \cite{W1}, is an important class of fully nonlinear elliptic operator. It is a generalization of Laplacian $(k=1)$ and Monge-Ampère $(k=n)$. 

First, for $k=1$,  (\ref{y101})-(\ref{y113}) become the classical Laplace equation
\begin{equation}\label{y107}
\begin{cases}
\Delta u=b(x) u^{\gamma}, & x \in \mathbb{R}^n, \\
u(x) \rightarrow \infty, & |x| \rightarrow \infty.
\end{cases}
\end{equation}
When $b(x) \equiv 1$,
  Keller \cite{K1} and Osserman \cite{O1}  obtained that the problem (\ref{y107}) has a positive radial solution $u \in C^2\left(\mathbb{R}^n\right)$ if and only if  the sublinear case  $0<\gamma< 1$.
When $b(x)\geq 0$ is a smooth function,
  Cheng and Ni \cite{CN1} studied the
  existence and uniqueness of the positive maximal entire  solution $U(x)$ for the superlinear case   $\gamma>1$, 
  and then established the asymptotic behavior 
  of 
  \begin{equation*}
  	U(x)\sim |x|^{(l-2)/(\gamma-1)},~|x|\rightarrow \infty,
  \end{equation*}
   if 
  $b(x) \sim|x|^{-l}$ with $l>2$ at $\infty$. Here, we define  $f \sim g$ at $\infty$  means that there exist positive constants $C_1$ and $C_2$, such that $C_1 f(x) \leq g(x) \leq C_2 f(x)$ for $|x|$ sufficiently large.
  In the sublinear case  $0<\gamma<1$, for radial locally Hölder continuous function $b(x)$ in $\mathbb{R}^n$,
  Lair and Wood \cite{L1}
  proved that (\ref{y107}) possesses an entire  radial  solution  if and only if
  \begin{equation*}
  \int_0^{\infty} r b(r) d r=\infty;
  \end{equation*}
  for continuous functions $b(x)$ which are not necessarily radial in $\mathbb{R}^n$,  assuming that
  \begin{equation}\label{y111}
  \int_0^{\infty} r b_{osc}(r) \exp \left(\int_0^r s b_*(s) d s\right) d r<\infty,
  \end{equation}
  where $b_{osc}(r):=b^*(r)-b_*(r)$, $b^*(r):=\sup _{|x|=r} b(x)$ and $b_*(r):=\inf_{|x|=r} b(x)$,  
  Lair \cite{L2} proved that (\ref{y107}) has an entire  solution if and only if
  \begin{equation}\label{y108}
  \int_0^{\infty} r b_*(r) d r=\infty;
  \end{equation}
  for locally bounded functions $b(x)$ which are not necessarily radial in $\mathbb{R}^n$, assuming that 
  \begin{equation}\label{y110}
  \int_0^{\infty} r b_{osc}(r) \left(1+\int_0^r s b_*(s) d s\right)^{\gamma/(1-\gamma)} d r<\infty,
  \end{equation}
  which is weaker than the condition (\ref{y111}),	
  El Mabrouk and Hansen \cite{EH} considered  that (\ref{y107}) has an entire  solution if and only if (\ref{y108}) holds, and then
  they also showed that (\ref{y107}) has no large solution if it admits a  positive bounded solution;
 Yang \cite{Y1} described the precise asymptotic behavior  of entire  solutions
 \begin{equation*}
 u(x) \sim\left\{\begin{array}{l}
 |x|^{(2-l)/(1-\gamma)}, ~ |x|\rightarrow\infty,~ \text { if } l<2, \\
 \left(\ln |x|\right)^{1/(1-\gamma)}, ~ |x|\rightarrow\infty,~ \text { if } l=2,
 \end{array}\right.
 \end{equation*}
 when  $b(x) \sim|x|^{-l}$ with $l \leq 2$ at $\infty$.

	Now let us return to (\ref{y101}),
	When $b(x) \equiv 1$ and $\gamma=0$,
 Bao, Chen, Guan and Ji \cite{BCGJ} proved that any convex solution $u \in C^{\infty}\left(\mathbb{R}^n\right)$ of (\ref{y101})	
	satisfying $u(x) \geq C_1|x|^2+C_2$   is a quadratic polynomial.
When $b(x) \equiv 1$ and  $\gamma>k$,	
 Jin, Li and Xu \cite{JHX} proved that   (\ref{y101}) has no positive entire $k$-convex  subsolution.  
 Later, Ji and Bao \cite{JB} showed that if $b(x) \equiv 1$, (\ref{y101}) has a positive entire  radial solution
 if and only if $0<\gamma< k$.
 
 More results on the existence, uniqueness and asymptotic behavior  of boundary blow-up solutions to (\ref{y101}) in the bounded domain  see
 \cite{CT,H1,J1,ML,S1,WSQ,ZF}. Other works of entire solutions to the Monge-Ampère  or Hessian equation, see \cite{CW,GT,JW,ZZ}.	
	The literatures for the sublinear case  $0<\gamma<k$   and  the entire  solution $u$ are less extensive.

In this paper, we  established the existence and exact asymptotic behavior of entire  $k$-convex solutions to the problem (\ref{y101})-(\ref{y113}). Our results are summarized as follows.

First, 
 we  assume that $b(x)$ is a positive continuous function in $\mathbb{R}^n$ and 
the oscillation function $b_{osc}(r)$ is small enough in the sense that
	\begin{equation}\label{y102}
	\begin{aligned}
\int_{0}^{\infty} \left(\frac{n r^{k-n}}{C_n^k} \int_{0}^{r}s^{n-1}b_{osc}(s) \tilde{b}(s) ds\right)^{1/k}
d r
<\infty,
\end{aligned}
\end{equation}
where 
\begin{equation*}
	\tilde{b}(s)=\left(1+\int_{0}^{s} \left(\frac{n t^{k-n}}{C_n^k} \int_{0}^{t} \tau^{n-1} b_*(\tau) d\tau \right)^{1/k} dt\right)^{k\gamma/(k-\gamma)},
\end{equation*}
 then we get a sufficient and necessary condition for the existence of entire   solution of the problem (\ref{y101})-(\ref{y113}). The proof of Theorem \ref{ythm1} is based on a nonexistence theorem for large solutions in a bounded domain and the sub-supersolution method.
\begin{theorem}\label{ythm1}
	Suppose that the condition (\ref{y102}) holds,
 then the problem (\ref{y101})-(\ref{y113}) admits positive entire  solutions $u\in\Phi^k(\mathbb{R}^n)$ if and only if
	\begin{equation}\label{y103}
	 \int_0^{\infty} \left(\frac{n r^{k-n}}{C_n^k} \int_{0}^{r} s^{n-1} b_*(s) ds\right)^{1/k} d r=\infty.
	\end{equation}
\end{theorem}
\begin{remark}
	the condition (\ref{y103}) is weaker than	
	\begin{equation}\label{y308}
	\int_{0}^{\infty} r  b_*^{1/k}(r) dr=\infty.
	\end{equation}	
	 For $k=1$, the condition (\ref{y103})  is equivalent to the condition (\ref{y108}), 
	   the condition (\ref{y102})  is equivalent to the condition (\ref{y110}). The detailed proof see Remark \ref{yrem301}.  Therefore, Theorem \ref{ythm1}  is   a generalization of  \cite{EH}.
\end{remark}
For the radial  function $b(x)$ in $\mathbb{R}^n$,  we find $b_{osc}(r)\equiv 0$. Then by Theorem \ref{ythm1}, we get the corollary below, which proves the existence of entire  radial solutions of the problem (\ref{y101})-(\ref{y113}).
\begin{corollary}\label{ythm2}
	Suppose that $b\left(x\right)=b\left(|x|\right) \in C\left(\mathbb{R}^n\right)$ is  positive, 
	then the problem (\ref{y101})-(\ref{y113}) admits positive entire  radial solutions  $u\in\Phi^k(\mathbb{R}^n)$  if and only if
	\begin{equation}\label{y109}
	\int_0^{\infty} \left(\frac{n r^{k-n}}{C_n^k} \int_{0}^{r} s^{n-1} b(s) ds\right)^{1/k} d r=\infty.
	\end{equation}
\end{corollary}
\begin{remark}
	For the sublinear case,
	if we extend $u^{\gamma}$ with $0<\gamma<k$ to $f(u) \in C^1(0, \infty)$ is a nonnegative, nondecreasing function, $f(0)=0$
	and
	$$\int_{1}^{\infty}f^{-1/k}(t) dt=\infty.$$
	Then we can generalize Theorem \ref{ythm1} and Corollary \ref{ythm2} and  get similar results (see \cite{Y1}).
\end{remark}
Finally, we establish  the exact asymptotic behavior of entire  solutions of the problem (\ref{y101})-(\ref{y113})  if $b(x)$ is a mild oscillation of a radial function $|x|^{-l}$.  Here, we define  $f=O(g)$ at $\infty$  means that there exist a positive constant $C$, such that $|f(x) / g(x)| \leq C$ for $|x|$ sufficiently large.
	\begin{theorem}\label{ythm3}
		  Suppose that  $b(x)=|x|^{-l}+O\left(|x|^{-m}\right)$ at $\infty$ for $l \leq k-1$, $m>l+(2k-l)k/(k-\gamma)$, then the problem (\ref{y101})-(\ref{y113}) admits  positive entire solutions $u\in\Phi^k(\mathbb{R}^n)$ satisfying
		 \begin{equation}\label{y104}
		 u(x)\sim |x|^{(2k-l)/(k-\gamma)},~  |x|\rightarrow \infty.
		 \end{equation}
		 Moreover, if $b(x)=b(|x|)$ is radial, then  $u(x)=u(|x|)$ is radial and
		 \begin{equation}\label{y105}
		 u^{\prime}(|x|)\sim |x|^{(2k-l)/(k-\gamma)-1},~  |x|\rightarrow \infty.
		 \end{equation}
		 \begin{equation}\label{y106}
		 u^{\prime \prime}(|x|) 
		 \sim |x|^{(2k-l)/(k-\gamma)-2},~  |x|\rightarrow \infty.
		 \end{equation}
	\end{theorem}
 The strategy for the proof of the asymptotic behavior comes from \cite{CN1}. The upper bound estimate is direct, and the lower bound estimate is inspired by the idea of \cite{NU}, which serves to control the behavior of the solutions at $\infty$.

The paper is organized as follows. In section 2,
 we give the comparison lemma that will be used in proving our results, then the proof of Corollary \ref{ythm2} and the nonexistence result for large solutions in a bounded domain  are given later. 
In section 3, we show a nonexistence result for large solutions if (\ref{y101}) admits a bounded solution and then  prove  Theorem \ref{ythm1}.
   In section 4, we give the proof of the asymptotic behavior for large solutions.

\section{Radial case}
In this section, we investigate  (\ref{y101}) in the case where $b(x)=b(|x|)$ is radially symmetric in $\mathbb{R}^n$. Our main goal here is to give a sufficient and necessary condition for the existence of an entire radial  solution to the problem (\ref{y101})-(\ref{y113})  (see Corollary \ref{ythm2}). 
This conclusion is useful for constructing  supersolutions and subsolutions  of (\ref{y101}) for the nonradial case in section 3.
The following comparison principle is a basic tool for proofs of the article and the detailed proof see \cite{U1}.
\begin{lemma}\cite{U1}\label{ylem302}
	(Comparison principle) 
	Suppose that   $u,v\in C(\bar{\Omega}) \cap \Phi^k(\Omega)$ are respectively the subsolution and supersolution of
	\begin{equation*}
		S_k\left(D^2 u\right)=\phi\left(x,u\right),
	\end{equation*}
	where $\phi\left(x,u\right)$ is positive
	and nondecreasing only with respect to $u$. If $u \leq v$ on $\partial \Omega$, then $u \leq v$ in $\Omega$.
\end{lemma}
\begin{corollary}\label{ycor301}
	Suppose that  $g\in C(\partial\Omega)$ is positive, $u,v\in C(\bar{\Omega}) \cap \Phi^k(\Omega)$ are	respectively the solutions of the following problems
	\begin{equation}\label{y304}
	\begin{cases}\sigma_k\left(\lambda\left(D^2 u\right)\right)=b(x) u^\gamma, & x\in  \Omega, \\ u=g, & x\in \partial \Omega\end{cases}
	\end{equation}
	and
	\begin{equation}\label{y305}
	\begin{cases}\sigma_k\left(\lambda\left(D^2 v\right)\right)=b(x) v^\gamma, & x\in \Omega, \\ v=\mu g, & x\in \partial \Omega,\end{cases}
	\end{equation}
	then $v\geq \mu u$ for $\mu\geq 1$.
\end{corollary}
\begin{proof}
	Since $\mu \geq 1$ and $\gamma<k$, we get from (\ref{y304}) that
	\begin{equation}\label{y306}
	\begin{cases}\sigma_k\left(\lambda\left(D^2 \left(\mu u\right)\right)\right)=\mu^k \sigma_k\left(\lambda\left(D^2 u\right)\right)
	\geq b(x) \left(\mu u\right)^\gamma, & x\in \Omega, \\ \mu u=\mu g, & x\in \partial \Omega.\end{cases}
	\end{equation}
	Using (\ref{y305}), (\ref{y306}) and Lemma \ref{ylem302}, we get $v \geq \mu u$ in  $\Omega$.
\end{proof}
To prove Corollary \ref{ythm2}, the nonexistence result for large solutions in a bounded domain  will be considered first.
\begin{lemma}\label{ylem201}
	Suppose that $\Omega\subseteq \mathbb{R}^n$ is a bounded domain and $b(x)$ 
	is a positive  continuous function on $\bar{\Omega}$, then the following problem 
		\begin{equation}\label{y207}
	\begin{cases}\sigma_k\left(\lambda\left(D^2 u\right)\right)=b(x) u^\gamma, &x\in \Omega, \\ u\rightarrow \infty, &x \rightarrow \partial \Omega\end{cases}
	\end{equation}
	has no positive  solution $u\in C(\bar{\Omega}) \cap \Phi^k(\Omega)$.
\end{lemma}
\begin{proof}
	On the contrary, suppose that (\ref{y207}) has a positive large solution $u(x)$. We choose $\alpha>0$, such that $g(x)=u(x)+\alpha\geq1$. Let
	\begin{equation}\label{y208}
	v(x)=\int_1^{u(x)+\alpha} s^{-\gamma/k} ds.
	\end{equation}
	By computation, we have
	\begin{equation*}
		v_i=
		\left(u+\alpha\right)^{-\gamma/k}
		 u_i,~
		v_{ij}
			=\left(u+\alpha\right)^{-\gamma/k}
			u_{ij}
			-\frac{\gamma}{k} (u+\alpha)^{-\gamma/k-1}
			u_i u_j,
	\end{equation*}
	i.e.
	\begin{equation*}
		D^2 v=\left(u+\alpha\right)^{-\gamma/k} D^2 u-\frac{\gamma}{k} (u+\alpha)^{-\gamma/k-1}	Du \otimes Du.
	\end{equation*}
		
	From \cite{ML}, we see, for any symmetric matrix $A$, $\xi \in \mathbb{R}^n$, $\lambda \in \mathbb{R}$, $1 \leq k \leq n$,
\begin{equation*}
S_k(A+\lambda \xi \otimes \xi)=S_k(A)+\lambda S_{k-1}\left(\left.A\right|_{i j}\right) \xi_i \xi_j, 
\end{equation*}
where $\left.A\right|_{i j}$ is the cofactor of the $(i,j)$-th entry of the symmetric matrix $A$.

 Therefore, considering $A=\left(u+\alpha\right)^{-\gamma/k} D^2 u$, $\lambda=-\frac{\gamma}{k} (u+\alpha)^{-\gamma/k-1}$, $\xi=Du$, we get
	\begin{equation*}
		 \begin{aligned}
		 S_k(D^2 v)&= S_k\left(
		 \left(u+\alpha\right)^{-\gamma/k}
		 D^2 u\right)-\frac{\gamma}{k} (u+\alpha)^{-\gamma/k-1} S_{k-1}\left(\left(u+\alpha\right)^{-\gamma/k} \left.D^2 u\right|_{i j}\right) u_i u_j\\
		 &=\left(u+\alpha\right)^{-\gamma} S_k(D^2 u)-\frac{\gamma}{k} (u+\alpha)^{-\gamma-1} S_{k-1}\left(\left.D^2 u\right|_{i j}\right) u_i u_j.
		  \end{aligned}
	\end{equation*}
Since $u$ is $k$-convex, then the matrix $D^2 u$ is positive definite. Therefore, we have $S_{k-1}\left(\left.D^2 u\right|_{i j}\right) u_i u_j\geq 0$ (see \cite{ZF}). It follows that
\begin{equation*}
	S_k(D^2 v)\leq \frac{S_k(D^2 u)}{\left(u+\alpha\right)^{\gamma}} = \frac{b(x) u^{\gamma}}{\left(u+\alpha\right)^{\gamma}}\leq b(x)\leq C,~x\in\Omega,
\end{equation*}	
	for some constant $C>0$. 
	
	Let $w(x)=M|x|^2\in C^2(\Omega)\cap C(\bar{\Omega})$, where $2M=\left(C_n^k\right)^{-1/k}C^{1/k}$. Then
	\begin{equation*}
		S_k(D^2 w)= C_n^k \left(2M\right)^k  \geq b(x).
	\end{equation*}
As a consequence, we have 
\begin{equation*}
	S_k(D^2 v)\leq S_k(D^2 w)~ \text {in}~ \Omega.
\end{equation*}	
Now let $\beta\geq 0$ be an arbitrary real number, then 
\begin{equation*}
S_k(D^2 v)\leq S_k\left(D^2 \left(w+\beta\right)\right) \text { in } \Omega.
\end{equation*}
By (\ref{y207}) and (\ref{y208}), $v(x) \rightarrow \infty$ as $x \rightarrow \partial \Omega$, it follows that
\begin{equation*}
	w+\beta\leq v,~ x \rightarrow \partial \Omega.
\end{equation*}
From the comparison principle \ref{ylem302}, we conclude that
\begin{equation*}
w+\beta\leq v,~ x \in \Omega,
\end{equation*}
for any real number $\beta\geq 0$. We choose $\beta$ sufficiently large and obtain a obvious contradiction.
\end{proof}

Next, considering  the case where $b(x)=b(|x|) \geq 0$ is radial in $\mathbb{R}^n$,  the existence of the positive entire radial solution of (\ref{y101})  is equivalent to the existence of the positive entire  solution of the Cauchy problem
	\begin{equation}\label{y201}
\begin{cases}
C_{n-1}^{k-1}u''\left(r\right) \left(\frac{u'\left(r\right)}{r}\right)^{k-1}
+C_{n-1}^{k}\left(\frac{u'\left(r\right)}{r}\right)^k=b(r)u^{\gamma}\left(r\right), ~r \in[0, \infty), \\
u^{\prime}(0)=0, ~ u(0)=a.
\end{cases}
\end{equation}
Integrating (\ref{y201}) from $0$ to $r$ yields
\begin{equation}\label{y209}	
u'(r)=\left(\frac{nr^{k-n}}{C_n^k} \int_0^r s^{n-1}  b(s) u^{\gamma}(s) d s\right)^{1 /k}. 
\end{equation}
From (\ref{y209}), we see that $u'\geq 0$. Therefore, $u$ is a nondecreasing function, and so $u(r)\geq a$ for $r\geq 0$.	
\begin{proof}[Proof of Corollary {\upshape\ref{ythm2}}]
First, we show that 	
 for any positive constant $a$, there exist a constant $R>0$, such that the Cauchy problem (\ref{y201}) has a solution $u(r)$ in $[0, R]$.
 
 We define a functional $F[\cdot, \cdot]$ on
 $$
 \mathcal{R}:=[0, R] \times\left\{u \in C[0, R]: a \leq u<2 a\right\},
 $$
 as
 \begin{equation*}
 F[r, u]:=\left(\frac{nr^{k-n}}{C_n^k} \int_0^r s^{n-1}  b(s) u^{\gamma}(s) d s\right)^{1 /k},
 \end{equation*}
 where $R$ is a small enough positive constant. Then (\ref{y201}) can be rewritten as
 \begin{equation*}
 u^{\prime}(r)=F[r, u] .
 \end{equation*}
 It is easy to see $F>0$ for $r>0$.
 
 We define an Euler's break line on $[0, R]$ as
 \begin{equation}\label{y202}
 \left\{\begin{array}{l}
 \psi(r)=a, ~0 \leq r \leq r_1, \\
 \psi(r)=\psi\left(r_{i-1}\right)+F\left[r_{i-1}, \psi\right]\left(r-r_{i-1}\right), ~ r_{i-1}<r \leq r_i,~ i=2,3, \cdots, m,
 \end{array}\right.
 \end{equation}
 where $0=r_0<r_1<\cdots<r_m=R$ and $m \in \mathbb{N}$.
 
 Step 1: We want to make sure that $a \leq \psi(r)<2 a$ for all $r \in[0, R]$, i.e. $(r,\psi)\in \mathcal{R}$.  
 In fact, it is obvious that $\psi(r) \geq a$. Since
 \begin{equation}\label{y203}
F\left[r_{i-1}, \psi\right]=\left(\frac{nr^{k-n}}{C_n^k} \int_0^r s^{n-1}  b(s) \psi^{\gamma}(s) d s\right)^{1 /k}
\leq
\frac{R  \max_{r\leq R} b^{1/k}(r) \psi^{\gamma/k}(R) }{\left(C_n^k\right)^{1/k}}.
 \end{equation}
Hence, for the break line $(r, \psi)$, we have
 \begin{equation*}
 a \leq \psi(r) 
 \leq a+\frac{R  \max_{r\leq R} b^{1/k}(r) \psi^{\gamma/k}(R) }{\left(C_n^k\right)^{1/k}} r 
 \leq 
 a+\frac{R^2  \max_{r\leq R} b^{1/k}(r) \psi^{\gamma/k}(R) }{\left(C_n^k\right)^{1/k}}.
 \end{equation*}
 Therefore, we can choose $R$ sufficiently small to make sure that $\psi(r)<2 a$.
 
 Step 2: We will prove that Euler's break line $\psi$ is an $\varepsilon$-approximation solution of (\ref{y201}). To do this, we only need to prove that for any small $\varepsilon>0$, there are appropriate points $\left\{r_i\right\}_{i=1, \cdots, m}$ to make the break line satisfy
 \begin{equation}\label{y204}
 \left|\frac{d \psi(r)}{d r}-F[r, \psi]\right|<\varepsilon, ~ r \in[0, R] .
 \end{equation}
 By (\ref{y202}), $\psi(r)$ has continuous derivatives in $[0, R]$ expect for a few points. There are unilateral derivatives at these individual points. If the derivative does not exist, we consider the right derivative.
 
 As a matter of fact, by (\ref{y203}), it is easy to see that
 \begin{equation*}
 \lim _{r \rightarrow 0} F[r, \psi]=0
 \end{equation*}
 is valid uniformly for any $(r, \psi) \in \mathcal{R}$. Then, for each $\varepsilon>0$, take $$\bar{r}=\frac{\left(C_n^k\right)^{1/k} \varepsilon}{\max_{r\leq R} b^{1/k}(r) \left(2a\right)^{\gamma/k}} \in(0, R),$$ 
  such that for $0 \leq r<\bar{r}$, we have
 $$
 F[r, \psi]<\varepsilon.
 $$
 Assume that $r_1=\bar{r}$, then
 \begin{equation*}
 \left|\frac{d \psi(r)}{d r}-F[r, \psi]\right|=|F[r, \psi]|<\varepsilon, ~ 0<r<\bar{r}.
 \end{equation*}
 For $\bar{r} \leq r \leq R$, let $r_{i-1}<r \leq r_i$, without loss of generality, we have
 \begin{equation*}
 \begin{aligned}
 & \left|\frac{d \psi(r)}{d r}-F[r, \psi]\right| \\
 \leq & \left|\frac{nr_{i-1}^{k-n}}{C_n^k} \int_0^{r_{i-1}} s^{n-1}  b(s) \psi^{\gamma}(s) d s
 -\frac{nr^{k-n}}{C_n^k} \int_0^r s^{n-1}  b(s) \psi^{\gamma}(s) d s\right|^{1 /k} \\
 \leq & \left(\frac{n\left|r_{i-1}^{k-n}-r^{k-n}\right|}{C_n^k} \int_0^{r_{i-1}} s^{n-1}  b(s) \psi^{\gamma}(s) d s
 -\frac{nr^{k-n}}{C_n^k} \int_{r_{i-1}}^r s^{n-1}  b(s) \psi^{\gamma}(s) d s\right)^{1 /k} \\
 \leq & \left(\frac{\left|r_{i-1}^{k-n}-r^{k-n}\right|}{C_n^k} R^{n}  \max_{r\leq R}b(r) (2a)^{\gamma} 
 -\frac{\bar{r}^{k-n}}{C_n^k} \left(r^n-r_{i-1}^n \right) \max_{r\leq R}b(r) (2a)^{\gamma}  \right)^{1 /k} .
 \end{aligned}
 \end{equation*}
 Since function $r^{k-n}$ and $r^n$ are both continuous on $[\bar{r}, R]$, for the above $\varepsilon$, we can take $\delta(\varepsilon)>0$ sufficiently small, such that
 $$
 \left|r^{\prime k-n}-r^{\prime \prime k-n}\right|<\left(2 R^{n}  \max_{r\leq R}b(r) (2a)^{\gamma} \right)^{-1} C_n^k \varepsilon^{k},
 $$
 $$
 \left|r^{\prime n}-r^{\prime \prime n}\right|
 <\left(2 \bar{r}^{k-n} \max_{r\leq R}b(r) (2a)^{\gamma}\right)^{-1} C_n^k \varepsilon^{k},
 $$
 where $r^{\prime}, r^{\prime \prime} \in[\bar{r}, R]$ and $\left|r^{\prime}-r^{\prime \prime}\right|<\delta(\varepsilon)$. 
 
 Assume that
 \begin{equation*}
 \max _{2 \leq i \leq m}\left|r_{i-1}-r_i\right|<\delta(\varepsilon),
 \end{equation*}
 then we get (\ref{y204}). Thus, Euler's break line $\psi$ is an $\varepsilon$-approximation solution of (\ref{y201}).

 Step 3: The next step is to find a solution of (\ref{y201}) by the Euler break line we defined. Assume  $\left\{\varepsilon_j\right\}_{j=1}^{\infty}$ is a positive constant sequence converging to $0$. For each $\varepsilon_j$, there is an $\varepsilon_j$-approximation solution $\psi_j$ on $[0, R]$, defined as above. By Step 1, it is easy to know that
  $$
 \left|\psi_j\left(r^{\prime}\right)-\psi_j\left(r^{\prime \prime}\right)\right|
 =F\left[r_{i-1}, \psi\right] \left|r^{\prime}-r^{\prime \prime}\right|
  \leq M\left|r^{\prime}-r^{\prime \prime}\right|,
 $$
 where $\left(r^{\prime},\psi_j\right),\left(r^{\prime\prime},\psi_j\right)  \in \mathcal{R}$. That is to say, $\left\{\psi_j\right\}$ is equicontinuous and uniformaly bounded $(r''=0)$. Then, by the Ascoli-Arzela Lemma, we can find a uniformly convergent subsequence, still denoted as $\left\{\psi_j\right\}$, without loss of generality.
 
 Assume that $\lim _{j \rightarrow \infty} \psi_j=u$. Since $\psi_j$ is continuous in $[0, R]$, we know that $u$ is continuous in $[0, R]$. By $\psi_j(0)=a$, we have $u(0)=a$.
 
 Since $\psi_j$ is an $\varepsilon_j$-approximation solution, we have
  \begin{equation}\label{y205}
 \frac{d \psi_j(r)}{d r}=F\left[r, \psi_j\right]+\Delta_j(r),
 \end{equation}
 where $\left|\Delta_j(r)\right|<\varepsilon_j$, for $r \in[0, R]$. Integrating (\ref{y205}) from $0$ to $r(\leq R)$, we have
 $$
 \psi_j(r)=a+\left(\int_0^r F\left[s, \psi_j\right] d s+\int_0^r \Delta_j(s) d s\right) .
 $$
 Let $j \rightarrow \infty$,
 \begin{equation}\label{y206}
 u(r)=a+\lim _{j \rightarrow \infty}\left(\int_0^r F\left[s, \psi_j\right] d s+\int_0^r \Delta_j(s) d s\right)=a+\int_0^r F[s, u] d s .
 \end{equation}
 Since $u$ is continuous in $[0, R]$, by (\ref{y206}), $u$ is continuously differentiable in $(0, R]$. Differentiating (\ref{y206}), we have $u^{\prime}(r)=F\left[r, u\right]$, $r>0$.
 Hence, we can see that $u$ satisfies (\ref{y201}) in $[0, R]$.

	 Next,  we show that 	
	 for any positive constant $a$,   (\ref{y201}) has a positive solution $u(r)$ in $[0, \infty)$.
	By (\ref{y209}), note that for $r \in\left[0, R\right]$, $0 \leq u^{\prime}(r)<\infty$.
	It implies that $u$ is nondecreasing and is finite if $r$ is finite. On the contrary,   suppose that $u$ blows up at some $r^*$, i.e. $u$ is a positive large  solution of (\ref{y101}) with $\Omega=B_{r^*}$, which
	contradicts to Lemma \ref{ylem201}. Therefore, $u$ can be extended to $[0, \infty)$ and remains positive. 
	
	Finally, we prove $u(r)\rightarrow \infty$ as $r\rightarrow \infty$ if and only if (\ref{y109}) holds. 	
	Assuming now that (\ref{y109}) does not hold, 
	We claim that $u$ is bounded.
	By (\ref{y209}), we know that		
	\begin{equation*}	
	u'(r)
	\leq
	\left(\frac{nr^{k-n}}{C_n^k} \int_0^r s^{n-1}  b(s)  d s\right)^{1 /k} u^{\gamma/k}(r)
	,
	\end{equation*}
	i.e. 
	\begin{equation}\label{y318}
		\frac{u'(r)}{u^{\gamma/k}(r)}\leq
		\left(\frac{nr^{k-n}}{C_n^k} \int_0^r s^{n-1}  b(s)  d s\right)^{1 /k}.
	\end{equation}
	Integrating (\ref{y318}) from $0$ to $r$, we get 
	\begin{equation*}
		u(r)\leq\left(
		a^{(k-\gamma)/k}+\frac{k-\gamma}{k}
		\int_0^{\infty} \left(\frac{n r^{k-n}}{C_n^k} \int_{0}^{r} s^{n-1} b(s) ds\right)^{1/k} d r
		\right)^{k/(k-\gamma)}<\infty.
	\end{equation*}	
	On the other hand, the integral form of (\ref{y201}) is
		\begin{equation*}
		\begin{aligned}	
		u(r)
		\geq 
		a+a^{\gamma/k} \int_0^r \left(\frac{ns^{k-n}}{C_n^k} \int_0^s t^{n-1}  b(t)  d t\right)^{1 /k} ds.
		\end{aligned} 
		\end{equation*}		
	Hence, combining (\ref{y109}), we get 	$u(r)\rightarrow \infty$ as $r\rightarrow \infty$. 		
\end{proof}	
\section{Nonradial case}
In this section, we consider the general case, where $b(x)$ is not necessarily radial. We present  conditions under which the problem (\ref{y101})-(\ref{y113}) has an entire  solution in $\mathbb{R}^n$. In order to prove Theorem \ref{ythm1}, we need the following lemmas.
First, we will discuss some sufficient conditions under which (\ref{y101}) cannot admit a positive large solution.
\begin{lemma}\label{ylem301}
 Suppose that $\Omega\in \mathbb{R}^n$ is a smooth domain and  
 $b(x)$ is a positive  continuous function in $\Omega$. Then (\ref{y101}) has no large solution in $\Omega$ if it has a  positive bounded solution in $\Omega$.
\end{lemma}
\begin{proof}
	 Let $u$ be a positive bounded solution of (\ref{y304}) and 
	 $$\max _{x \in \bar{\Omega}} u(x)=c_0>0.$$
	 Choose $c \geq c_0$ and then $\mu:=\frac{c}{c_0}\geq 1$.  Moreover, by Corollary \ref{ycor301}, we obtain that for any $c\geq c_0$,
	 \begin{equation*}
	 	v\geq \frac{c}{c_0} u,
	 \end{equation*}
	 where $v$ is the solution  of (\ref{y305}) and
	 the proof of the existence is similar to Theorem 3.1 in \cite{BDS}. 	 
Suppose now that (\ref{y101}) has a large solution $w$ in $\Omega$. Then, by the comparison principle in Lemma \ref{ylem302}, $w\geq v$ and thereby $w \geq \frac{c}{c_0} u$ in $\Omega$ for any $c \geq c_0$. This yields a contradiction because $u$ is not identically zero in $\Omega$.
\end{proof}
Next, we give the comparison principle of the initial  value.
	\begin{lemma}\label{ylem402}
	Suppose that $u_1(r)$ is a positive entire solution of		
	\begin{equation*}
	\begin{cases}
	C_{n-1}^{k-1}u_1''\left(\frac{u_1'}{r}\right)^{k-1}+C_{n-1}^{k}\left(\frac{u_1'}{r}\right)^k=b_1(r)u_1^{\gamma}, ~r \in[0, \infty), \\
	u_1^{\prime}(0)=0, ~ u_1(0)=\alpha_1>0
	\end{cases}
	\end{equation*}
	and $u_2(r)$ is a positive entire solution of
	\begin{equation*}
	\begin{cases}
	C_{n-1}^{k-1}u_2''\left(\frac{u_2'}{r}\right)^{k-1}+C_{n-1}^{k}\left(\frac{u_2'}{r}\right)^k=b_2(r)u_2^{\gamma}, ~r \in[0, \infty), \\
	u_2^{\prime}(0)=0, ~ u_2(0)=\alpha_2 \geq
	 \alpha_1,
	\end{cases}
	\end{equation*}
	where  $b_2(r) \geq b_1(r) \geq 0$ in $[0, \infty)$. Then $u_2(r) \geq u_1(r)$ in $[0, \infty)$.
\end{lemma}
\begin{proof}
	Suppose by contradiction that the conclusion is false, then there exists $r_0>0$, such that  $u_2(r_0)=u_1(r_0)$ and $u_2'(r_0)<u_1'(r_0)$. This implies that  $u_2(r)\geq u_1(r)$ in $[0,r_0]$. 
	By direct calculation, we have	
	\begin{equation*}
	\begin{aligned}
	u_2'(r_0)&= \left(\frac{nr^{k-n}}{C_n^k} \int_0^{r_0} s^{n-1} b_2(s)  u_2^{\gamma}(s) d s\right)^{1 /k}\\
	&\geq
	\left(\frac{nr^{k-n}}{C_n^k} \int_0^{r_0} s^{n-1} b_1(s)  u_1^{\gamma}(s) d s\right)^{1 /k}
	=u_1'(r_0),
	\end{aligned}
	\end{equation*}
	we obtain a contradiction.
	The proof of this lemma is completed.
\end{proof}
By Lemma \ref{ylem402}, we can estimate the radial solution of (\ref{y101})  for $b(x)=b(r)$, $r=|x|$.
	\begin{lemma}\label{ylem303}
	Suppose that  $u(r)$ is the positive entire radial  solution of  
	$$\sigma_k\left(\lambda\left(D^2 u\right)\right)
	=  b(r) u^\gamma,$$
	then
	\begin{equation}\label{y307}
	u \leq 2^{\gamma/(k-\gamma)}\left(u(0)+
	\overline{u}
	^{k/(k-\gamma)}\right),
	\end{equation}
	where $\overline{u}$ is the positive entire radial solution of $\sigma_k\left(\lambda\left(D^2 \overline{u}\right)\right)
	=  b(r)$. 
\end{lemma}
\begin{proof}
	By direct calculation, we have
	\begin{equation*}
		\overline{u}(r)=\int_0^r\left( \frac{ns^{k-n}}{C_n^k} \int_0^s t^{n-1}  b(t) d t \right)^{1/k}d s.
	\end{equation*}
	Let $\varepsilon>0$, $f:=(u(0)+\varepsilon)^{(k-\gamma)/k}+\overline{u}$ and $\tilde{u}:=f^{k/(k-\gamma)}$, then we get
		\begin{equation*}
	f'=\overline{u}',~	f''=\overline{u}'',
	\end{equation*}
	\begin{equation*}
	\tilde{u}'=\frac{k}{k-\gamma}f^{\gamma/(k-\gamma)}f',~
	\tilde{u}''=\frac{k}{k-\gamma}f^{\gamma/(k-\gamma)}f''+
	\frac{k\gamma}{(k-\gamma)^2}f^{(2\gamma-k)/(k-\gamma)} (f')^2,
	\end{equation*}
	and
	\begin{equation*}
	\begin{aligned}
	 \sigma_k\left(\lambda\left(D^2 \tilde{u}\right)\right)&=C_{n-1}^{k-1}\tilde{u}''\left(\frac{\tilde{u}'}{r}\right)^{k-1}+C_{n-1}^{k}\left(\frac{\tilde{u}'}{r}\right)^k\\
	 &=
	 \frac{k^k f^{k\gamma/(k-\gamma)}}{\left(k-\gamma\right)^k}	 
	 \left(C_{n-1}^{k-1}f''\left(\frac{f'}{r}\right)^{k-1}
	 +C_{n-1}^{k}\left(\frac{f'}{r}\right)^k
	 +\frac{C_{n-1}^{k-1} \gamma (f')^{k+1}}{(k-\gamma)fr^{k-1}}\right)\\
	 &=
	 	 \frac{k^k \tilde{u}^{\gamma}}{\left(k-\gamma\right)^k}		 	 
	 \left( b(r)
	 +\frac{C_{n-1}^{k-1} \gamma (f')^{k+1}}{(k-\gamma)fr^{k-1}}\right)
	  \\
	  &\geq  b(r) \tilde{u}^{\gamma}.
	 \end{aligned}
	\end{equation*}	
	So, since $u(0) \leq \tilde{u}(0)$ and $\tilde{u}>0$ in $\mathbb{R}^n$, $u(x) \leq \tilde{u}(x)$ for all $x \in \mathbb{R}^n$ by Lemma \ref{ylem402}. Thus
	\begin{equation*}
	u \leq\left((u(0)+\varepsilon)^{(k-\gamma)/k}+\overline{u}\right)^{k/(k-\gamma)} \leq 2^{\gamma/(k-\gamma)}\left(u(0)+\varepsilon+\overline{u}^{k/(k-\gamma)}\right).
	\end{equation*}
	Letting $\varepsilon\rightarrow0$, we get (\ref{y307}).
\end{proof}
\begin{remark}
	When $b(r)$ is replaced by $b_*(r)$ or $b^*(r)$, Lemma \ref{ylem303} still holds.
\end{remark}
	\begin{proof}[Proof of Theorem {\upshape\ref{ythm1}}]
		First, we will prove the sufficient condition.
 By the proof of Corollary \ref{ythm2}, we know that for any $\beta>1$, the following problems
 		\begin{equation}\label{y301}
 \begin{cases}
 C_{n-1}^{k-1}v''\left(\frac{v'}{r}\right)^{k-1}+C_{n-1}^{k}\left(\frac{v'}{r}\right)^k=b^*(r)v^{\gamma}, ~r \in[0, \infty), \\
 v^{\prime}(0)=0, ~ v(0)=1
 \end{cases}
 \end{equation}
 and
 		\begin{equation}\label{y302}
\begin{cases}
C_{n-1}^{k-1}w''\left(\frac{w'}{r}\right)^{k-1}+C_{n-1}^{k}\left(\frac{w'}{r}\right)^k=b_*(r)w^{\gamma}, ~r \in[0, \infty), \\
w^{\prime}(0)=0, ~ w(0)=\beta
\end{cases}
\end{equation}
have positive  solutions $v(r)$ and $w(r)$ in $[0, \infty)$. Obviously, $v$ is a subsolution and $w$ is a supersolution of  (\ref{y101}). According to Theorem 2.10 of \cite{N1}, (\ref{y101}) has a positive solution $u$, such that $v \leq u \leq w$ in $[0, \infty)$ if 
\begin{equation}\label{y303}
v \leq w, ~ r \in[0, \infty)
\end{equation}
holds for $\beta$ sufficiently large.

To show (\ref{y303}), we first estimate $w$. Set $\underline{w}=\int_1^w t^{-\gamma / k} d t$, then $w=\left(1+\frac{k-\gamma}{k} \underline{w}\right)^{k /(k-\gamma)}$. Similar to Lemma \ref{ylem303}, we have
\begin{equation*}
\underline{w}^{\prime}=w^{-\gamma / k} w^{\prime},~ 
\underline{w}^{\prime \prime}=w^{-\gamma / k} w^{\prime \prime}-\frac{\gamma}{k} w^{-\gamma / k-1}\left(w^{\prime}\right)^2 \leq w^{-\gamma / k} w^{\prime \prime} .
\end{equation*}
Then $\underline{w}$ satisfies
\begin{small}
\begin{equation}\label{y309}
\begin{aligned}
C_{n-1}^{k-1} \underline{w}^{\prime \prime}\left(\frac{\underline{w}^{\prime}}{r}\right)^{k-1}
\!+\!
C_{n-1}^k\left(\frac{\underline{w}^{\prime}}{r}\right)^k  \!\leq\!
 w^{-\gamma}\left(C_{n-1}^{k-1} w^{\prime \prime}\left(\frac{w^{\prime}}{r}\right)^{k-1}
 \!+\!
 C_{n-1}^k\left(\frac{w^{\prime}}{r}\right)^k\right) \!=\!
 b_*(r) .
\end{aligned}
\end{equation}
\end{small}
Integrating (\ref{y309}) from $0$ to $r$, we infer
\begin{equation*}
\underline{w}(r) \leq \int_0^r\left(\frac{n s^{k-n}}{C_n^k} \int_0^s t^{n-1}  b_*(t) d t\right)^{1 / k} d s.
\end{equation*}
Therefore, we get
\begin{equation}\label{y314}
w(r) \leq\left(1+\frac{k-\gamma}{k}
\int_0^r\left(\frac{n s^{k-n}}{C_n^k} \int_0^s t^{n-1}  b_*(t) d t\right)^{1 / k} d s
\right)^{k /(k-\gamma)}.
\end{equation}
Now, suppose by contradiction that (\ref{y303}) is not true. Then there exists a constant $r_*>0$, such that $v(r)<w(r)$ for $r \in\left(0, r_*\right)$ and $v\left(r_*\right)=w\left(r_*\right)$. Integrating (\ref{y301}) and (\ref{y302}) from 0 to $r_*$, we have
\begin{equation}\label{y311}	
v(r_*)=1+\int_{0}^{r_*} \left(\frac{ns^{k-n}}{C_n^k} \int_0^s t^{n-1}  b^*(t) v^{\gamma}(t) d t\right)^{1 /k} d s
\end{equation}
and
\begin{equation}\label{y312}	
w(r_*)=\beta+\int_{0}^{r_*} \left(\frac{ns^{k-n}}{C_n^k} \int_0^s t^{n-1}  b_*(t) w^{\gamma}(t) d t\right)^{1 /k} d s.
\end{equation}
By  (\ref{y311}) and (\ref{y314}), we infer
\begin{small}
\begin{equation}\label{y310}
\begin{aligned}	
&v(r_*)\leq 1
+\int_{0}^{r_*} \left(\frac{ns^{k-n}}{C_n^k} \int_0^s t^{n-1}  b_*(t) w^{\gamma}(t) d t\right)^{1 /k} d s +\\
&
\int_{0}^{r_*} \left(\frac{ns^{k-n}}{C_n^k} \int_0^s t^{n-1}  b^*(t) w^{\gamma}(t) d t\right)^{1 /k} d s
- \int_{0}^{r_*} \left(\frac{ns^{k-n}}{C_n^k} \int_0^s t^{n-1}  b_*(t) w^{\gamma}(t) d t\right)^{1 /k} d s\\
&\leq 1
+\int_{0}^{r_*} \left(\frac{ns^{k-n}}{C_n^k} \int_0^s t^{n-1}  b_*(t) w^{\gamma}(t) d t\right)^{1 /k} d s\\
&+\int_{0}^{\infty} \left(\frac{nr^{k-n}}{C_n^k} \int_0^r s^{n-1} b_{osc}(s) w^{\gamma}(s) d s\right)^{1 /k} 
d r
,
\end{aligned}
\end{equation}	
\end{small}
where 
\begin{equation*}
w^{\gamma}(s) \leq\left(1+
\int_0^s\left(\frac{n t^{k-n}}{C_n^k} \int_0^t \tau^{n-1}  b_*(\tau) d \tau\right)^{1 / k} d t
\right)^{k\gamma /(k-\gamma)}=\tilde{b}(s).
\end{equation*}
By condition  (\ref{y102}), it is easy to see 
 (\ref{y310}) is finite.
 Thus $\beta$ may be chosen sufficiently large, such that $v\left(r_*\right)<w\left(r_*\right)$. A contradiction is obtained.  Using (\ref{y303}), the proof of the sufficient condition will be completed if we show that $\lim _{|x| \rightarrow \infty} u(x)=\infty$. 
By (\ref{y303}), (\ref{y311}) and (\ref{y103}), we  see that 
\begin{equation}\label{y316}
\begin{aligned}
u(x) \geq v(r) \geq 1+ \int_{0}^{r} \left(\frac{ns^{k-n}}{C_n^k} \int_0^s t^{n-1}  b_*(t)  d t\right)^{1 /k} d s 
\rightarrow \infty, ~|x|\rightarrow \infty.
\end{aligned}
\end{equation}

Second, we will prove the necessary condition by contradiction.
Suppose now that (\ref{y103}) does not hold, we claim that $u$ is bounded. By (\ref{y303}), (\ref{y312}) and using Lemma \ref{ylem303}, we get
	\begin{equation*}
u \leq w \leq 2^{\gamma/(k-\gamma)}\left(\beta+
\overline{w}
^{k/(k-\gamma)}\right),
\end{equation*}
where
\begin{equation}\label{y313}
\begin{aligned}
\overline{w}(r) =\int_0^r\left( \frac{ns^{k-n}}{C_n^k} \int_0^s t^{n-1}  b_*(t) d t \right)^{1/k}d s.
\end{aligned}
\end{equation}
since (\ref{y103}) does not hold, then (\ref{y313}) trivially implies that
 $w$ is bounded which yields that $u$ is also bounded. Since $u$ is a nontrivial entire bounded  solution, then from Lemma \ref{ylem301}, (\ref{y101}) has no  entire large solution.
	\end{proof}
\begin{remark}
	From the proof of Theorem \ref{ythm1}, we have that	
		assuming  the condition (\ref{y102}) holds,
	then (\ref{y101}) admits positive entire bounded solutions $u\in\Phi^k(\mathbb{R}^n)$ if and only if
		\begin{equation*}
	\int_0^{\infty} \left(\frac{n r^{k-n}}{C_n^k} \int_{0}^{r} s^{n-1} b_*(s) ds\right)^{1/k} d r<\infty.
	\end{equation*}
 On the other hand, if $b$ satisfies the condition (\ref{y103}), then (\ref{y101}) has no positive entire bounded solution in $\mathbb{R}^n$.
\end{remark}
	\begin{remark}\label{yrem301}
		First, we  show that the condition (\ref{y103}) is weaker than	the condition (\ref{y308}).	
	 Using the fact that the function $\phi(t)= t^{1/k}$ is concave, then applying Jensen's inequality, we obtain
	\begin{equation*}
	\left(\frac{\int_{0}^{s} t^{n-1} b_*(t) d t}{\int_{0}^{s} t^{n-1} dt}\right)^{1/k}
	\geq 
	\frac{\int_{0}^{s} t^{n-1} b_*^{1/k}(t) d t}{\int_{0}^{s} t^{n-1} dt},
	\end{equation*}
	then 
	\begin{equation*}
	\left(n s^{k-n} \int_{0}^{s} t^{n-1} b_*(t) d t\right)^{1/k}\geq 
	n s^{1-n} \int_{0}^{s} t^{n-1} b_*^{1/k}(t) d t.
	\end{equation*}
	By (\ref{y316}) and using Fubini's theorem, the solutions $u$ of (\ref{y101}) satisfy
	\begin{equation}\label{y317}
	\begin{aligned}	
	\liminf_{|x|\rightarrow \infty}u(x)
	&\geq 1+ 	\frac{1}{\left(C_n^k\right)^{1/k}}
	\int_{0}^{\infty} \left(ns^{k-n} \int_0^s t^{n-1}  b_*(t)  d t\right)^{1 /k} d s \\
	&\geq 1+
	\frac{n}{\left(C_n^k\right)^{1/k}}
	\int_{0}^{\infty} 
	s^{1-n} \int_0^s t^{n-1}  b_*^{1/k}(t)  d t d s\\
	&\geq 1+	
	\frac{n}{\left(C_n^k\right)^{1/k}}
	\int_{0}^{\infty} 
	t^{n-1} b_*^{1/k}(t)  \int_t^r  s^{1-n}  d s d t\\
	&\geq 1+	\frac{n}{\left(n-2\right)\left(C_n^k\right)^{1/k}}
	\int_{0}^{\infty} 
	t b_*^{1/k}(t)   d t.
	\end{aligned}
	\end{equation}
	Hence, 
	if the condition (\ref{y308}) holds, we can deduce that the condition (\ref{y103}) holds.
 In particular, for $k=1$, the condition (\ref{y103})  is equivalent to the condition (\ref{y108}).
 Notice that 
 \begin{equation}\label{y315}
 	\int_0^r s^{1-n} \int_0^s t^{n-1} b_*(t) dt ds
 	=\frac{1}{n-2} \left(\int_0^r sb_*(s)ds-r^{2-n}\int_0^r s^{n-1}b_*(s) ds \right).
 \end{equation}
Then by  (\ref{y303}), Lemma \ref{ylem303} and (\ref{y315}), we have	
\begin{equation}\label{y319}
	\begin{aligned}	
	u(x)\leq w(r)&\leq 2^{\gamma/(1-\gamma)}
	\left(\beta+\left(\int_0^r s^{1-n} \int_{0}^s t^{n-1} b_*(t) dtds\right)^{1/(1-\gamma)}\right)\\
	&\leq 2^{\gamma/(1-\gamma)}
	\left(\beta+\left(\frac{1}{n-2}\int_0^r s b_*(s) ds\right)^{1/(1-\gamma)}\right).
		\end{aligned}	
\end{equation}
Conditions  (\ref{y317}) and (\ref{y319})	
 guarantee that conditions (\ref{y103}) and (\ref{y108}) are equivalent.

	Next, we  show that the condition (\ref{y102}) is stronger than 
\begin{equation}\label{y114}
\int_0^{\infty} r b_{osc}^{1/k}(r) \left(1+\frac{n}{\left(C_n^k\right)^{1/k}}\int_0^r s^{1-n} \int_0^s t^{n-1} b_*^{1/k}(t) dt ds\right)^{\gamma/(k-\gamma)} dr<\infty.
\end{equation}		
Similar to (\ref{y317}), applying Jensen's inequality repeatedly, we obtain
 \begin{equation*}
 	\int_{0}^{\infty} \left(\frac{n r^{k-n}}{C_n^k} \int_{0}^{r}s^{n-1}b_{osc}(s) \tilde{b}(s) ds\right)^{1/k}
 	d r
 	\geq \frac{n}{(n-2)\left(C_n^k\right)^{1/k}} 
 	\int_0^{\infty} r b_{osc}^{1/k}(r) \tilde{b}^{1/k}(r) dr,
 \end{equation*}
where
\begin{equation*}
	\begin{aligned}
	\tilde{b}^{1/k}(r)
	&\geq \left(1+\frac{n}{\left(C_n^k\right)^{1/k}}\int_0^r s^{1-n} \int_0^s t^{n-1} b_*^{1/k}(t) dt ds\right)^{\gamma/(k-\gamma)}.
	\end{aligned}
\end{equation*}
	Hence, 
if the condition (\ref{y102}) holds, we can deduce that the condition (\ref{y114}) holds.
In particular, for $k=1$, since the condition (\ref{y103})  is equivalent to the condition (\ref{y108}), we can find the condition (\ref{y102})  is equivalent to the condition (\ref{y110}) easily  and the proof is similar to  \cite{L1}.	
	\end{remark}
\begin{remark}
	Notice that (\ref{y102}) is not necessary. For example, for $k=1$,  $b(x)=8\left(2 x_1^2+x_2^2+x_3^2+1\right)^{-1/2}$ and $\gamma=1/2$,  the problem (\ref{y101})-(\ref{y113}) has the nonradial solution $u(x)=2 x_1^2+x_2^2+x_3^2+1$ and $b_{osc}(r)=O(1 / r)$ as $r \rightarrow \infty$, such that  (\ref{y102}) does not hold (see \cite{L2}).
\end{remark}	
	\section{Asymptotic behavior}
	In this section, we study the asymptotic behavior of solutions for the problem (\ref{y101})-(\ref{y113}) and complete the proof of Theorem \ref{ythm3}. We seperate it into several steps, and the radial case as follows will be considered first.
	\begin{lemma}\label{ylem401}
		Suppose that $u(r)$, $r=|x|$ is a positive radial  solution of the problem (\ref{y101})-(\ref{y113}) with $b(x)=b(r)$  in $\mathbb{R}^n$ and $b(r)=r^{-l}$ in $\mathbb{R}^n \backslash B_{R_0}$ for some $R_0>0$ and $l \leq k-1$, 		
	then  
	  $$u(r)\sim r^{(2k-l)/(k-\gamma)},~ r\rightarrow \infty.$$
\end{lemma}
\begin{proof}
	First, we prove that there exists a constant $C>0$, such that  $u(r) \leq C r^{(2k-l)/(k-\gamma)}$,  $r\rightarrow \infty$.	
	By (\ref{y201}), $u(r)$ satisfies 
	\begin{small}
				\begin{equation}\label{y401}
		\frac{C_n^k r^{1-n}}{n}\left(r^{n-k}\left(u'\left(r\right)\right)^k\right)'
		= C_{n-1}^{k-1}u''\left(r\right) \left(\frac{u'\left(r\right)}{r}\right)^{k-1}
		+ C_{n-1}^{k}\left(\frac{u'\left(r\right)}{r}\right)^k= b(r)u^{\gamma}\left(r\right).
		\end{equation}
	\end{small}
Integrating (\ref{y401}) from $0$ to $r$ yields
	\begin{equation}\label{y402}
		r^{n-k} \left(u'\left(r\right)\right)^k=  \int_0^r \frac{n s^{n-1}}{C_n^k}  b(s)  u^{\gamma}(s) d s.
	\end{equation}
	So $u(r)$ is nondecreasing for $r\geq 0$, and then
	\begin{equation}\label{y403}
		\int_0^r \frac{n s^{n-1}}{C_n^k}  b(s)  u^{\gamma}(s) d s
		\geq 
		u^{\gamma}(0) \int_{R_0}^r \frac{n s^{n-1-l}}{C_n^k}   d s \rightarrow \infty, ~ r \rightarrow \infty.
	\end{equation}
	By (\ref{y402}) and (\ref{y403}), we have
		\begin{equation*}
	\begin{aligned}
	r^{n-k} \left(u'\left(r\right)\right)^k&=  \int_0^{R_0} \frac{n s^{n-1}}{C_n^k}  b(s)  u^{\gamma}(s) d s+
	\int_{R_0}^r \frac{n s^{n-1-l}}{C_n^k}    u^{\gamma}(s) d s\\
	&\leq C_1
	\int_{R_0}^r \frac{n s^{n-1-l}}{C_n^k}    u^{\gamma}(s) d s\\
 &\leq
 \frac{n C_1}{(n-l) C_n^k} r^{n-l} u^{\gamma}(r)
 ,
	\end{aligned}
	\end{equation*}
	i.e.
	\begin{equation}\label{y404}
			\frac{u^{\prime}(r)}{u^{\gamma/k} (r)} \leq  \left(\frac{n C_1}{\left(n-l\right) C_n^k}\right)^{1/k} r^{1-l/k},
	\end{equation}
	where $C_1>1$. Here, we choose  $R_0>0$ is small, such that $C_1$ is close to $1$.	
 Integrating (\ref{y404}) from $R_0$ to $r$, we get that for $r$ sufficiently large, 
	\begin{equation}\label{y405}
	u(r)  \leq C_2 r^{(2k-l)/(k-\gamma)}.
	\end{equation}
Next, we prove that there exists a constant $C>0$, such that  $u(r) \geq C r^{(2k-l)/(k-\gamma)}$,  $r\rightarrow \infty$.	
By 	(\ref{y401}) and (\ref{y404}), we have for $r>R_0$  sufficiently large,
	\begin{equation}\label{y408}
		C_{n-1}^{k-1}u''(r)\left(\frac{u'(r)}{r}\right)^{k-1}\geq
		C_3 r^{-l} u^{\gamma}(r),
	\end{equation}
	where $C_3=1-C_1(n-k)/(n-l)\in \left(0, 1\right]$.
	
And then by  (\ref{y408}) and (\ref{y405}), we have
\begin{equation}\label{y406}
u''(r)\left(u'(r)\right)^{k-1}
\geq C_3 \left(C_{n-1}^{k-1}\right)^{-1} r^{k-l-1} u^{\gamma}\left(r\right)
\geq
	 C_4 u^{(k-\gamma)(k-l-1)/(2k-l)+\gamma}(r).
\end{equation}
Multiplying by $u'(r)$ on (\ref{y406}) and integrating from $R_0$ to $r$, we get that for $r$  sufficiently large,

\begin{equation*}
	\left(u'(r)\right)^{k+1}\geq C_5 u^{(k-\gamma)(k-l-1)/(2k-l)+\gamma+1}(r),
\end{equation*}
i.e.
\begin{equation}\label{y407}
	u'(r)\geq C_5 u^{1-(k-\gamma)/(2k-l)}(r).
\end{equation}
 Integrating (\ref{y407}) from $R_0$ to $r$ again, we get that for $r$  sufficiently large, 
 \begin{equation}\label{y409}
 u(r)\geq C_6 r^{(2k-l)/(k-\gamma)}.
 \end{equation}
 The proof  of this lemma is completed.
\end{proof}
\begin{remark}
	Under conditions of Lemma \ref{ylem401},
	 by (\ref{y405}) and (\ref{y409}), we know that the problem (\ref{y101})-(\ref{y113}) has a radial solution
	\begin{equation}\label{y410}
	u(x)=C |x|^{(2k-l)/(k-\gamma)},~x\in\mathbb{R}^n \backslash B_{R_0}.
	\end{equation}
	Since (\ref{y410}) satisfies (\ref{y401}), by direct calculation,  we have
	\begin{equation*}
	C=\left(\frac{n}{C_n^k\left(n+\left(\alpha-2\right)k\right)\alpha^k}\right)^{1/(k-\gamma)}>0,~\text{where}~\alpha=\frac{2k-l}{k-\gamma}.
	\end{equation*}
\end{remark}

	Now, we conclude the asymptotic behavior of radial solutions for the problem (\ref{y101})-(\ref{y113}) in the radial case  as follows.
	\begin{proposition}\label{yprop401}
		 Suppose that $u(r)$ is a positive radial  solution of the problem (\ref{y101})-(\ref{y113}) with  $b(x)=b(r)$ in $\mathbb{R}^n$ and $b(r)\sim r^{-l}$ at $\infty$ for $l \leq k-1$. Then (\ref{y104})-(\ref{y106}) hold for $u(r)$.
	\end{proposition}
	\begin{proof}
		Since $b(r)\sim r^{-l}$,   $r\rightarrow \infty$, there exist two radial functions $b_i(r)$, $i=1,2$, such that $b_1(r) \leq b(r) \leq b_2(r)$ in $[0, \infty)$ and $b_i(r)=C_i r^{-l}$ in some $\left[r_0, \infty\right)$, where $r_0>0$ and $C_i>0$, $i=1,2$.  By Corollary \ref{ythm2}, the equation 
		$$\sigma_k\left(\lambda\left(D^2 u\right)\right)= b_i(r) u^\gamma,~x\in\mathbb{R}^n$$ 
		 has a positive radial large  solution $u_i(r)$, $i=1,2$ with $u_1(0)= u(0)/2$ and $u_2(0)=2 u(0)$. Therefore, using Lemma \ref{ylem402} and Lemma \ref{ylem401} shows that $u_1(r) \leq u(r) \leq u_2(r)$ in $[0, \infty)$, and
		\begin{equation*}
		u_i(r)\sim r^{(2k-l)/(k-\gamma)},~i=1,2,~ r\rightarrow \infty.
		\end{equation*}
 This implies that (\ref{y104}) holds for $u(r)$. It remains to estimate $u^{\prime}(r)$ and $u^{\prime \prime}(r)$. Integrating (\ref{y401}) from $0$ to $r$, we infer
 	\begin{equation*}
 u'(r)\sim \left(r^{k-n} \int_{r_0}^r s^{n-1} s^{-l}  s^{\gamma(2k-l)/(k-\gamma)} d s\right)^{1 /k},~r\rightarrow\infty,
 \end{equation*}
 i.e.
 $$u'(r)\sim r^{(2k-l)/(k-\gamma)-1},~r\rightarrow\infty.$$
		Furthermore, by (\ref{y401}),
		\begin{equation*}
		\begin{aligned}
			u''(r)&=\left(C_{n-1}^{k-1}\right)^{-1}\left(b(r)r^{k-1}\left(u'(r)\right)^{1-k} u^{\gamma}(r)-C_{n-1}^k r^{-1}u'(r)\right)\\
			&\sim r^{(2k-l)/(k-\gamma)-2},~
			r\rightarrow\infty.
			\end{aligned}
		\end{equation*}
		The proof of this proposition is completed.
	\end{proof}
	\begin{proof}[Proof of Theorem {\upshape\ref{ythm3}}]
		 	According to the conditions of Theorem \ref{ythm3}, suppose that $b(x)=|x|^{-l}+O\left(|x|^{-m}\right)$ at $\infty$ for $l \leq k-1$, $m>l+(2k-l)k/(k-\gamma)$, then
		 we have $ b_*(r)\sim r^{-l}$ and $b_{osc}(r)=O\left(r^{-m}\right) $ at $\infty$. Obviously (\ref{y103}) holds.
		 
		 On the other hand, it is easy to verify that
		 \begin{equation*}
		 \left(r^{k-n} \int_0^r s^{n-1} b_{osc}(s) \tilde{b}(s) d s\right)^{1 /k}  
		 \leq C r^{ \left(k-m+\left(2k-l\right) \gamma/\left(k-\gamma\right) \right)/k} 
		 \leq C r^{-1},~r\rightarrow \infty,
		 \end{equation*}
		 then (\ref{y102})   holds. Therefore, by Theorem \ref{ythm1}, the problem (\ref{y101})-(\ref{y113}) has  positive entire  solutions.	 
		 Moreover, $v(|x|) \leq u(x) \leq w(|x|)$ in $[0, \infty)$, where $v$ satisfies (\ref{y301}) and $w$ satisfies (\ref{y302}) (see the proof of Theorem \ref{ythm1}). It is easy to see that $b^*(r)$ and $b_*(r) \sim r^{-l}$ at $\infty$ for $l \leq k-1$. By using Proposition \ref{yprop401}, (\ref{y104})-(\ref{y106}) hold for $v$, $w$, and then for $u$. The proof of Theorem \ref{ythm3} is completed.
		\end{proof}

\end{document}